\documentclass[letterpaper, 10 pt, conference]{ieeeconf}  

\IEEEoverridecommandlockouts                              

\usepackage[utf8]{inputenc}
\usepackage{amsmath,amssymb,amsfonts}
\def\disp{\displaystyle}

\def\oa{\bar a}

\def\ox{\bar{x}}

\def\ou{\bar{u}}

\def\l{\left}
\def\r{\right}

\def\gg{\gamma}

\def\nn{\left \{ }

\def\la{\langle}
\def\ra{\rangle}

\def\lm{\lambda}

\def\sr{\Rightarrow }

\def\th{\theta}
\def\vth{\vartheta}
\def\ox{\bar x}
\def\oa{\bar a}
\def\vt{\overrightarrow}

\newcommand{\R}{\mathbb{R}}

\newtheorem{theorem}{Theorem}[section]

\usepackage{graphicx,color}

\include{Xlatin1}

\newcommand{\ba}{\begin{array}}
\newcommand{\ea}{\end{array}}

\title{\LARGE \bf Applications of Controlled Sweeping Processes to Nonlinear Crowd Motion Models with Obstacles}

\author{Tan H. Cao$^{1}$, Boris S. Mordukhovich$^{2}$, Dao Nguyen$^{2}$, Trang Nguyen$^{2}$
\thanks{$^{1}$T.H. Cao is with Department of Applied Mathematics and Statistics, SUNY (State University of New York) Korea
    {\tt\small tan.cao@stonybrook.edu}}%
    \thanks{$^{2}$B.S. Mordukhovich, D. Nguyen, and T. Nguyen are with Department of Mathematics, Wayne State University, USA
        {\tt\small boris@math.wayne.edu, dao.nguyen2@wayne.edu, daitrang.nguyen@wayne.edu}}%
}

\begin{document}
\maketitle
\thispagestyle{empty}
\pagestyle{empty}
\begin{abstract}
This paper mainly focuses on solving the dynamic optimization of the planar controlled crowd motion models with obstacles which is an
application of a class of optimal control problems governed by a general perturbed nonconvex sweeping process. This can be considered as a significant extension of the previous work regarding the controlled crowd motion models in \cite{cm2,cm1}, where the obstacles have not been considered. The necessary optimality conditions for the problem under consideration are established based on the results in \cite{cm2,cm1} and are illustrated by a nontrivial example of practical importance.
\end{abstract}
\normalsize

\section{Introduction and Problem Formulation}\label{sec:intro}
This paper is mainly devoted to applications of a class of optimal control problems governed by a perturbed nonconvex sweeping process in \cite{cm1} to solving dynamic optimization problems involving the crowd motion models. More precisely, we consider the motion of $N$ agents without overlapping and aim at finding an optimal strategy to steer them to a desired target with the minimum effort during a fixed interval time $[0,T]$. Following the mathematical framework in \cite{mv2}, we identify $N$ agents to inelastic disks with different radii $R_i$ whose centers are denoted by $x_i = (x_{i1}, x_{i2})\in \R^2$ for $i=1,\ldots, N$. Our purpose is to simulate the motion of these agents inside a given domain $\Omega \subset \R^2$ during the given interval time $[0,T]$ and subject to both static and dynamic obstacles. To ensure the nonoverlapping of agents and obstacles, the {\it configuration vector} $x=(x_1,x_2,\ldots,x_N)\in \R^{2N}$ has to belong to the following set of {\it admissible configurations}
\begin{equation}
    \label{nov}
    C = \l\{x\in \R^{2N}|\; D_{ij}(x)\ge 0, \;\; \forall i< j,\; i, j\in\{1,\ldots, N\}\r\}
\end{equation}
where $D_{ij}(x):= \l\| x_i-x_j\r\| - (R_i+R_j)$ is the signed distance between disks $i$ and $j$. In this setting, each of $N$ agent has his/her own {\it desired velocity} in case they are not in contact in the sense that $D_{ij}>0$. When they are in contact, i.e., $D_{ij}=R_i+R_j$; however, the agents $i$ and $j$ must adjust their velocities in order to avoid overlapping. In this very situation the actual velocities and desired velocities are no longer the same while they should be close to each other. To reflect this situation, we introduce the set of {\it admissible velocities} $V(x)$ defined by
\begin{equation}
    \label{av}
    \l\{ v\in \R^{2N}|\; D_{ij}(x)=0 \sr \la G_{ij}(x),v\ra\ge 0,\; \forall i<j\r\},
\end{equation}
where
$$
\begin{aligned}
G_{ij}(x)&:=\nabla D_{ij}(x) \\
&=\l(0,\ldots, 0, -e_{ij}(x), 0,\ldots, 0, e_{ij}(x), 0\ldots, 0 \r)
\end{aligned}
$$
and where
$$
e_{ij}:=\dfrac{x_j-x_i}{\l\| x_j-x_i\r\|}.
$$
Next, let $U(x):=(v(x_1),v(x_2), \ldots, v(x_N))$ be the (desired) velocity of all the agents assuming that they behave the same behavior and reach the targets using the shortest path. In this case $v(x_i)=-s_i\nabla D(x_i)$, where $D(x_i):= \l\|x_i - x^{tar}_i \r\|$ stands for the distance between the agent $x_i$ and the corresponding $x_i$ and where $s_i$ denotes the speed.

To ensure the nonoverlapping of all the agents and obstacles, their actual velocity at time $t$ denoted by $\dot x(t)$ should be selected from the set of admissible velocities $V(x(t))$ in \eqref{av}. Let us see how this selection guarantees the nonoverlapping condition. Using the first-order Taylor approximation of the distance function $D_{ij}(\cdot)$ gives us
$$
D_{ij}(x(t+h)) = D_{ij}(x(t)) + h \la\nabla D_{ij}(x(t)),\dot x(t) \ra + o(h)
$$
where $h>0$ is a small unit of time. Suppose that agents $i$ and $j$ are in contact at time $t$, then $D_{ij}(x(t))=0$. Moreover, since $\dot x(t)\in V(x(t))$, then
$$
\la\nabla D_{ij}(x(t)),\dot x(t) \ra \ge 0,
$$
which implies that  $D_{ij}(x(t+h))\ge 0$. Therefore, the nonoverlapping condition in \eqref{nov} does not allow the agents to move with their desired velocity freely and the set of admissible velocities in \eqref{av} offers a guideline to select the actual velocity in order to avoid the nonoverlapping situation. The next step is how to select a ``good" velocity over the set of admissible velocities in \eqref{av}. In fact the actual velocity of $N$ agents is defined to be the closet one to their desired velocity in the following sense
\begin{equation}
    \label{act-vel}
    \dot x(t) = \Pi (U(x(t)); V(x(t)))
\end{equation}
where $\Pi (U(x(t)); V(x(t))):=\inf_{v\in V(x(t))}\l\| v-U(x(t))\r\|$ denotes the Euclidean projection of $U(x(t))$ onto the set $V(x(t))$. The differential equation in \eqref{act-vel} can be written in the equivalent form of perturbed sweeping process as follows using the orthogonal decomposition via the sum of mutually polar cone as in \cite{mv2}:
\begin{equation}
    \label{sp}
    \dot x(t) \in -N_C(x(t)) + U(x(t))
\end{equation}
where $N_C(x)$ stands for the {\it proximal normal cone} of the set $C$ at $x$. In fact, the sweeping process (``processes de raffle") was first introduced by Jean-Jacques Moreau in the 1970s in the differential form
\begin{equation}
    \label{msp}
    \dot x(t) \in - N_{C(t)}(x(t))\;\; \mbox{a.e.}\;\; t\in [0,T]
\end{equation}
where $C(t)$ is a continuous time moving convex set and the symbol $N_{\Omega}(x)$ signifies the classical normal cone to a convex set $\Omega$ at $x$ defined by
$$
N_{\Omega}(x):= \l\{
\begin{array}{llll}
    \l\{ v\in\R^n|\; \la v,u-x\ra \;\; \forall u\in \Omega \r\}& \mbox{if} & x\in \Omega\\
    \emptyset &\mbox{if} &x\notin \Omega
\end{array}
\r.
$$
The sweeping process \eqref{msp} describes the motion of an object inside the moving set $C(t)$. Depending on the behavior of the set $C(t)$, the object will stay where it is (in the case, where it is not hit by the boundary of $C(t)$), or otherwise it is swept towards the interior of the set $C(t)$. In the latter case, the velocity of the object points inwards to the moving set in order not to leave. The sweeping (Moreau) process was originally motivated by applications to elastoplasticity \cite{mor_frict}, but it has been further significantly recognized in the much broader spectrum of applications such as hysteresis, electric circuits, traffic equilibria, populations motion in confined spaces, and other areas of applied sciences and operations research. Its well-posedness was established by using the {\it catching-up} algorithm developed by Moreau \cite{mor_frict}. Since then, there have been several improvements of \eqref{msp} by many authors including, e.g., relaxing the convexity of the moving set $C(t)$, and extending \eqref{msp} to the perturbed version
\begin{equation}
 \label{psp}
 \dot x(t)\in  -N_{C(t)}(x(t)) + f(x(t)) \;\;\mbox{a.e.}\;\; t\in [0,T],
\end{equation}
where $f(\cdot)$ represents some given external force. In particular, Maury and Venel \cite{mv2} developed a mathematical framework for uncontrolled microscopic model of the crowd motion dynamic by using an extended sweeping process in form of \eqref{psp} and providing  its numerical simulations with various applications. Let us emphasize that the differential inclusion \eqref{sp} in the crowd motion framework is a special case of the perturbed sweeping process \eqref{psp} without involving control actions and optimization.

Motivated by control-theoretic requirements and practical applications, the authors of \cite{cm1} considered the following optimal control problem:
\begin{equation}
   \label{(P)}
   \mbox{minimize} \;\; J[x,a]: = \frac{1}{2}\l\| x(T) \r\|^2 +\frac{1}{2}\int^T_0\l\| a(t)\r\|^2\; dt
\end{equation}
over the control functions $a(\cdot)\in L^2([0,T];\R^N)$ and the corresponding trajectories $x(\cdot)\in W^{1,2}([0,T];\R^{2N})$ of the controlled nonconvex sweeping process
\begin{equation}
   \label{csp}
   \begin{cases}
   \dot x(t) \in -N_C(x(t)) + U(x(t),a(t))\;\; \mbox{a.e.}\;\; t\in [0,T],\\
   x(0) = x_0\in C,
   \end{cases}
\end{equation}
 where the set $C$ is given in \eqref{nov}, where the controlled desired velocity $U(x(t),a(t))$ is given by
 \begin{equation}
 \label{dv}
 U(x(t),a(t)) = (-a_1s_1\nabla D(x_1), \ldots, -a_Ns_N\nabla D(x_N)),
 \end{equation}

The meaning of the cost function in \eqref{(P)} is minimizing the distance of all the agents to the exit situated at the origin together with the minimum energy of feasible controls $a(\cdot)$ used to adjust the desired velocity. The authors of \cite{cm1,cm2} discussed how to solve such an optimal control problem for the planar crowd motion models with $N$ agents in general by using the obtained necessary optimality conditions. Observe that these efficient necessary conditions allow us to compute the optimal solutions analytically and successfully. However, no obstacles (either static or dynamic) have been considered in the controlled crowd motion model of \cite{cm1}. This paper is devoted to the study of the dynamic optimization of the controlled planar crowd motion models {\em subject to obstacles}. In this setting, all the agents consider neither overlapping with each others nor with the obstacles.

To reflect this situation, the sweeping set $C$ should be replaced by the following set $C^{obs}$ involving the obstacles:
\begin{equation}
\label{nov-obs}
\begin{aligned}
C^{obs}:=  &\big\{x\in \R^{2N}|\; D_{ij}(x)\ge 0, \;\; \forall i< j,\; i, j\in \{1,\ldots, N\} \\
 & D^{obs}_{ij}(x)\ge0\;\; \forall i\in\{1,\ldots, N\}, j\in\{1, \ldots, m\} \big\},
\end{aligned}
\end{equation}
where $D^{obs}_{ij}(x):= \l\| x_i - x^{obs}_j\r\| - (R_i + r_j)$ for $i \in\{1,\ldots, N\}$ and $j\in\{1, \ldots, m\} $.  Here, the $j$th obstacle is represented by a rigid disk in $\R^2$ with radius $r_j$ whose center is denoted by $x^{obs}_j$. Note that the sweeping set $C$ in \eqref{nov} cannot be used to represent the nonoverlapping conditions between the agents and the obstacles, since--if an agent is close enough to an obstacle he/she will avoid it--while the obstacle has no reaction to the agent. The obstacle always stays where it is or moves in its own direction without interacting with the agent. Therefore, the obstacles cannot be involved in the sweeping process representing the dynamics of all the agents in \eqref{csp}, and hence the sweeping set $C^{obs}$ in our framework is significantly different from the set $C$ in \eqref{nov}. In this paper we consider the following optimal control problem:
\begin{equation}
   \label{(P1)}
   \mbox{minimize} \;\; J[x,a]: = \frac{1}{2}\l\| x(T) - x^{tar} \r\|^2 +\frac{\tau}{2}\int^T_0\l\| a(t)\r\|^2\; dt
\end{equation}
over the control functions $a(\cdot)\in L^2([0,T];\R^N)$ and the corresponding trajectories $x(\cdot)\in W^{1,2}([0,T];\R^{2N})$ of the controlled nonconvex sweeping process
\begin{equation}
   \label{csp}
   \begin{cases}
   \dot x(t) \in -N_{C^{obs}}(x(t)) + U(x(t),a(t))\;\; \mbox{a.e.}\;\; t\in [0,T],\\
   x(0) = x_0\in C^{obs},
   \end{cases}
\end{equation}
 where the set $C^{obs}$ is given in \eqref{nov-obs}, where the controlled desired velocity $U(x(t),a(t))$ is given by \eqref{dv}, where $x^{tar}\in\R^{2N}$ stands for the desired target that we want to send all the agents to, and where $\tau\ge 0$ is a given constant.
 The objective of all the agents here is to minimize their distances to the desired targets by using the least energy. The scalar $\tau\ge0$ introduced in \eqref{(P1)} is a trade-off between the distance and the energy.

The paper is organized as follows. Section~II contains a set of necessary optimality conditions for problem $(P)$. In Section~III we provide numerical implementation of the obtained results for problem $(P)$ for a practical situation. The last Section~IV concerns topics of our future research.

We refer the reader to \cite{ac,cmn1,cmn,pfs,zeidan} and the bibliographies therein for  further recent developments on controlled sweeping processes and therein applications.\vspace*{-0.1in}

\section{Crowd Motion Models with Obstacles}
In this section we summarize the necessary optimality conditions to the optimization of the planar crowd motion model developed in paper \cite{cm1}.
\begin{theorem}
\label{Th1}
{\bf (necessary optimality conditions for optimization of controlled crowd motions)}
    Let $(\ox(\cdot), \oa(\cdot)) \in W^{2,\infty}([0,T];\R^{2N}\times \R^N)$ be a strong local minimizer of the crowd motion problem in \eqref{csp}. It follows from \cite[Theorem~3.1]{cm1} that there exist $\lm \ge 0, \eta_{ij}(\cdot) \in L^2([0,T];\R_+)$, for $i, j=1, \ldots, N$, well-defined at $t=T$, $w(\cdot) = (w^x(\cdot),w^a(\cdot)) \in L^2([0,T];\R^{2N}\times\R^N)$, $v(\cdot) = (v^x(\cdot),v^a(\cdot))\in L^2([0,T];\R^{2N}\times\R^N) $, a measure $\gamma \in C^*([0,T];\R^{2N})$, and a vector function $q(\cdot) = (q^x(\cdot),q^a(\cdot)):[0,T]\to \R^{2N}\times \R^N$ of bounded variation on $[0,T]$ such that the following conditions are satisfied:
    \begin{enumerate}
        \item[{\bf(1)}] $w(t) = (0, \oa(t)),\;\; v(t)=(0,0)$ for a.e. $t\in [0,T]$;

        \item[{\bf(2)}]
$\dot\ox(t)+(a_1s_1\nabla D(x_1(t)), \ldots, a_Ns_N\nabla D(x_N(t)))$\\[1ex]
$=\dot\ox(t) + (a_1s_1d(\ox_1(t),x^{tar}_1), \ldots, a_Ns_Nd(\ox_N(t),x^{tar}_N))$\\[1ex]

$\begin{aligned}
=&\disp\bigg(-\sum_{j>1}\eta_{1j}(t) d(\ox_j(t),\ox_1(t)),\ldots, \\
&\sum_{i<j}\eta_{ij}(t)d(\ox_j(t),\ox_i(t))
\disp-\sum_{i>j}\eta_{ji}(t)d(\ox_i(t),\ox_j(t)),\\
&\disp\ldots,\sum_{j<N}\eta_{jN}(t)d(\ox_N(t),\ox_j(t))\bigg),
\end{aligned}$\\[1ex]
where $d(x,y):=\dfrac{x-y}{\l\| x-y\r\|}$;\\[1ex ]

\item[{\bf(3)}] $\|\ox_i(t)-\ox_j(t)\|>2R\Longrightarrow\eta_{ij}(t)=0$ for all $i<j$ and a.e.\ $t\in[0,T]$;
\item[\bf{(4)}] $\eta_{ij}(t)>0\Longrightarrow\left\la q^x_j(t)-q^x_i(t),\ox_j(t)-\ox_i(t)\right\ra=0$ for all $i<j$ and a.e.\ $t\in[0,T]$;
\item[{\bf(5)}]
$\left\{\begin{array}{ll}
\dot p(t)=\bigg(0,\lm\oa_1(t)-s_1\big[\la q^x_1(t),d(\ox_1(t),x^{tar}_1)\ra \big],\\ \ldots,\lm\oa_N(t)-s_N\big[\la q^x_N(t),d(\ox_N(t),x^{tar}_N)\ra\big]
\bigg);\end{array}\right.$ \\[1ex]
\item[{\bf(6)}] $q^x(t)=p^x(t)+\gg([t,T])$ for a.e.\ $t\in[0,T]$;
\item[{\bf(7)}] $q^a(t)=p^a(t)=0$ for a.e.\ $t\in[0,T]$;
\item[{\bf(8)}]$\left\{\begin{array}{ll}
p^x(T)+\lm\ox(T)=\disp\bigg(-\sum_{j>1}\eta_{1j}(T)d(\ox_j(T),\ox_1(T)),\\[3ex]
\ldots,\disp\sum_{i<j}\eta_{ij}(T) d(\ox_j(T),\ox_i(T))\\[3ex]
\disp-\sum_{i>j}\eta_{ji}(T)d(\ox_i(T),\ox_j(T)),\ldots,\\[3ex]
\disp\sum_{j<N}\eta_{jN}(T)d(\ox_N(T), \ox_j(T)) \bigg);
\end{array}\right.\\[1ex]
$\item[\bf{(9)}] $p^a(T)=0$;
\item[{\bf(10)}] $\lm+\|p^x(T)\|>0$.
\end{enumerate}

\end{theorem}\vspace*{0.05in}

Rewriting the differential equations in (2) gives us
\begin{equation}\label{dyn1}
\l\{
\begin{array}{ll}
    \dot \ox_1(t) = -s_1\oa_1(t)d(\ox_1(t),x^{tar}_1)\\
     -\disp \sum_{j>1}\eta_{1j}(t)d(\ox_j(t),\ox_1(t))\\[1ex]
    \dot \ox_i(t) = -s_i\oa_i(t)d(\ox_i(t), x^{tar}_i)\\
     +\disp\sum_{i<j}\eta_{ij}(t)d(\ox_j(t),\ox_i(t)) \\[1ex]
    - \disp\sum_{i>j}\eta_{ji}(t)d(\ox_i(t),\ox_j(t)), i=2, \ldots, N-1,\\
    \dot \ox_N(t) = -s_N\oa_N(t)d(\ox_N(t),x^{tar}_N) \\
    + \disp\sum_{j<N}\eta_{jN}(t)d(\ox_N(t),\ox_j(t))
\end{array}
\r.
\end{equation}
for a.e. $t\in [0,T]$. In fact, this system of differential equations tells us that the actual velocity of each agent depends on the agent's current location with respect to the other agents and the target.

It follows from (5) and (7) that
\begin{equation}
\label{control1}
\lambda \oa_i(t) = s_i\la q^x_i(t), d(\ox_i(t), x^{tar}_i)\ra
\end{equation}
for $i=1,\ldots, N$. This equation gives us useful information about the controls relating to the trajectories. The crowd motion model with $N=2$ agents was fully studied in \cite{cm1}, where several nontrivial examples were considered in various settings and solved analytically. Our next goal is to consider the crowd motion models with obstacles.

\section{Numerical Implementations}
\begin{figure}[ht]
\centering
\includegraphics[scale=0.38]{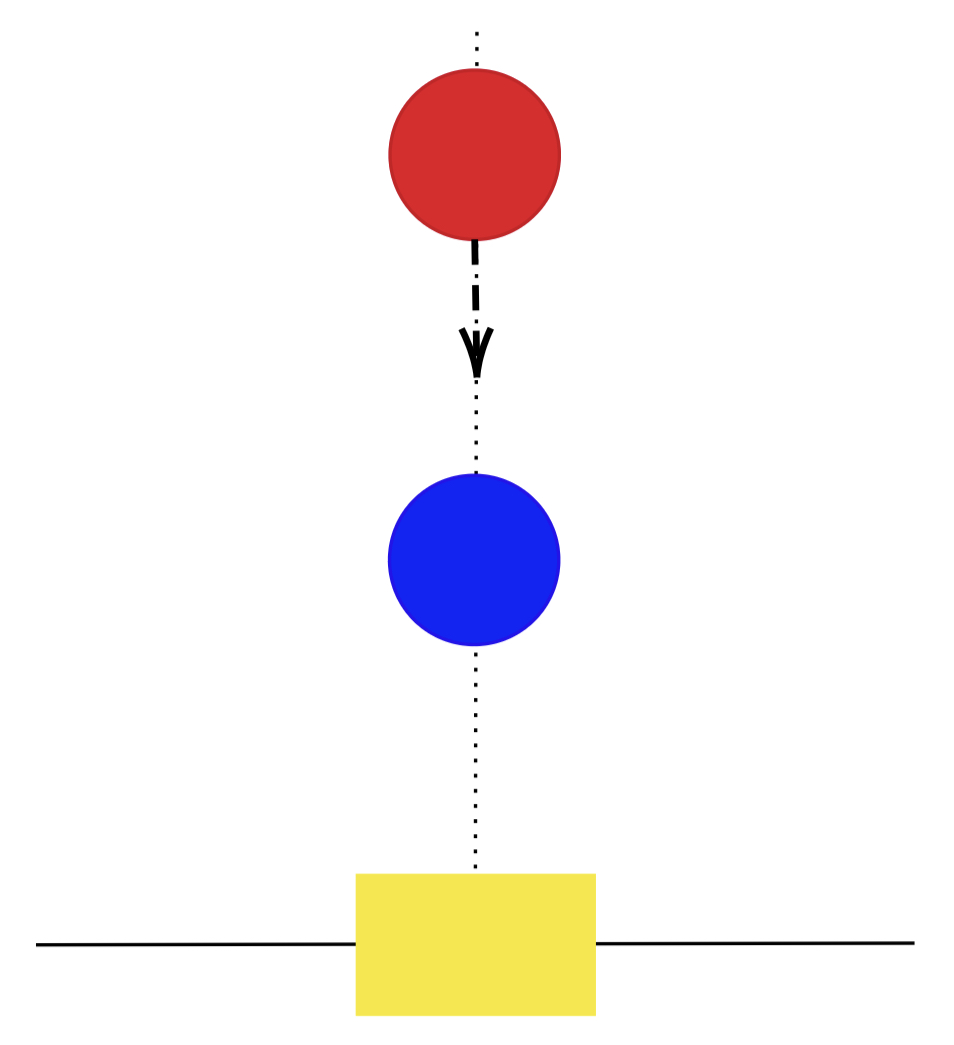}
\caption{The agent is moving towards the target but he/she must avoid the obstacle.}
\label{pic1}
\end{figure}

This section is devoted to the study of the crowd motion models with static obstacles. For simplicity, we first consider a single agent whose objective is to reach the target while avoiding the obstacle using the minimum energy. In this case, the obstacle is simply a state constraint involved in the sweeping set $C^{obs}$ defined by
$$
C^{obs}:=\{x\in R|\; D_1(x):=\l\| x-x^{obs}\r\| - 2R\ge 0\},
$$
where $x^{obs}$ represents the position of the obstacle in $\R^2$. Modifying the proof of necessary optimality conditions in Theorem~\ref{Th1} we can find some dual elements $\lm\ge 0, \eta_1(\cdot)\in L^2([0,T];\R_+)$ well-defined at $t=T$
$w(\cdot) = (w^x(\cdot),w^a(\cdot)) \in L^2([0,T];\R^{2}\times\R)$, $v(\cdot) = (v^x(\cdot),v^a(\cdot))\in L^2([0,T];\R^{2}\times\R) $, a measure $\gamma \in C^*([0,T];\R^{2})$, and a vector function $q(\cdot) = (q^x(\cdot),q^a(\cdot)):[0,T]\to \R^{2}\times \R$ of bounded variation on $[0,T]$ such that the following conditions are satisfied:

\begin{enumerate}
    \item[{\bf(1)}] $w(t) = (0, \oa(t)),\;\; v(t)=(0,0)$ for a.e. $t\in [0,T]$;
    \item[{\bf(2)}] $\dot\ox(t) = -s\oa(t)d\l(\ox(t),x^{tar}\r)-\eta(t)d\l(x^{obs},\ox(t)\r)$
    for a.e. $t\in [0,T]$;
    \item[{\bf(3)}] $\l\| \ox(t)-x^{obs}\r\|>2R\Longrightarrow\eta(t)=0$;
    \item[{\bf(4)}] $\eta(t)>0\Longrightarrow\left\la q^x(t),\ox^{obs} -\ox(t)\right\ra=0$ for a.e.\ $t\in[0,T]$;
    \item[{\bf(5)}]
    $\dot p(t)=\bigg(0,\lm\oa(t)-s\big[\la q^x(t),d(\ox(t),x^{tar})\ra \big]$;
        \item[{\bf(6)}] $q^x(t)=p^x(t)+\gg([t,T])$ for a.e.\ $t\in[0,T]$;
        \item[{\bf(7)}] $q^a(t)=p^a(t)=0$ for a.e.\ $t\in[0,T]$;
        \item[{\bf(8)}] $p^x(T)+\lm\ox(T)=-\eta(T)d(\ox^{obs},\ox(T))$;
        \item[\bf{(9)}] $p^a(T)=0$;
        \item[{\bf(10)}] $\lm+\|p^x(T)\|>0$.
\end{enumerate}
In what follows, we consider an example illustrating our necessary conditions obtained above. Assume that the target is located at the origin i.e, $x^{tar} =(0,0)$, the initial position of the agent is $(0,48)$ and the location of the static obstacle is $(0,24)$, see Figure~\ref{pic1}. The agent aims to reach the target while avoiding the obstacle using the minimum energy. The data of the problem can be described as follows:
$$
\l\{
\begin{array}{ll}
    T = 6, \;x_0 = (0, 48), \; x^{obs} = (0, 24), \; ,x^{tar}=(0,0); \\[1ex]
    s = \frac{48}{6} = 8, \;R = 3,\; \ell = \frac{1}{2}\l(\oa^2 \r) , \varphi(x) = \frac{1}{2}x^2.
\end{array}
\right.
$$
It follows from condition {\bf(2)} that
$$
\dot\ox(t) = -8\oa(t)\dfrac{\ox(t)}{\l\| \ox(t) \r\|} - \eta(t)\dfrac{x^{obs}-\ox(t)}{\l\|x^{obs}-\ox(t) \r\|}
$$
for all $t\in [0,T]$. Let $t_{1}$ and $\vartheta_{1}$ be the first time and last time that the agent and the obstacle are in contact, that is,
$$\|x^{obs}(t)-\ox(t)\|=2R=6$$
for all $t\in [t_{1},\vartheta_{1}]$. Then the velocity of the agent is given by
\begin{equation}\label{dynamic}
 \dot\ox(t) = \l\{
  \begin{array}{llll}
   -\dfrac{8\oa(t)}{\l\|\ox(t)\r\|} (\ox_1(t), \ox_2(t)), &\mbox{if } t \notin [t_1,\vartheta_1]\,\\[2ex]
   -\dfrac{8\oa(t)}{\l\|\ox(t)\r\|} (\ox_1(t), \ox_2(t)) \\
    -\dfrac{\eta(t)}{6}(-\ox_1(t),24-\ox_2(t)), &\mbox{if } t\in [t_1,\vartheta_1],\\[2ex]
  \end{array}
  \r.
\end{equation}
which is based on conditions {\bf(2)} and {\bf(3)}. Before being in contact with the obstacle, the agent moves towards the target using his/her own desired velocity. However, when he/she is in contact with the obstacle, the normal cone $N_{C^{obs}}(\ox(t))$ in \eqref{csp} is activated and hence generates a normal vector $-\eta(t)\dfrac{(x^{obs}-\ox(t))}{\l\|x^{obs}-\ox(t) \r\|}$ forcing the agent to move away from the obstacle.
Let us investigate behavior of the agent when he/she is in contact with the obstacle during the time interval $[t_1,\vartheta_1]$. If $t\in [t_1,\vartheta_1]$, then
$\|x^{obs}(t)-\ox(t)\|=2R=6$, which implies that
$$\la\ox(t)-x^{obs},\ox(t)-x^{obs}\ra=4R^2.$$
Differentiating both sides of the above equation with respect to $t$ gives us
\begin{equation}
\label{ort}
\la\dot{\ox}(t),\ox(t)-x^{obs}\ra =0.
\end{equation}
Combining it with the representation of $\dot{\ox}(t)$ in \eqref{dynamic} we obtain
$$
\bigg\la-8\oa(t)\dfrac{\ox(t)}{\|\ox(t)\|}-\eta(t)\dfrac{x^{obs}-\ox(t)}{\|x^{obs}-\ox(t)\|},\;\ox(t)-x^{obs}\bigg\ra =0,
$$
and hence get the very useful information about the scalar function $\eta(\cdot)$ as follows
\begin{equation}\label{eta}
\begin{aligned}
\eta(t) &=\dfrac{\oa(t)}{4R\ \|\ox(t)\|}\la\ox(t),\ox(t)-x^{obs}\ra \\
&= \frac{1}{2}\oa(t)\bigg\la \dfrac{\ox(t)}{\|\ox(t)\|}, \dfrac{\ox(t)-x^{obs}}{\|\ox(t)-x^{obs}\|} \bigg\ra.
\end{aligned}
\end{equation}
This shows that the behavior of $\eta(\cdot)$ depends heavily on the position of the agent with respect to the target and the obstacle. When the agent first meets the obstacle, two unit vectors $\dfrac{\ox(t_1)}{\|\ox(t_1)\|}$ and $\dfrac{\ox(t_1)-x^{obs}}{\|\ox(t_1)-x^{obs}\|}$ are parallel which causes the maximum effect to push the agent away from the obstacle. In this case, we have $\eta(t_1) = \dfrac{1}{2}\oa(t_1)>0$. In fact, due to structure of the problem, the control function $\oa(\cdot)$ takes the positive values on $[0,T]$. For simplicity, we assume that $\oa(\cdot)$ is constant on $[0,T]$, i.e., $\oa(t) = \oa$ for all $t\in[0,T]$. Once the agent is in contact with the obstacle, he/she would stay on the circle centered at $x^{obs}=(0,24)$ with radius $2R=6$ during the time interval $[t_1,\vartheta_1]$. Hence, it makes a perfect sense to assume that $\eta(t)>0$ on $[t_1,\vartheta_1)$.

We then deduce from the implication in ${\bf(4)}$ that $q^x(t)$ is orthogonal to $\ox(t)-x^{obs}$, which implies that $q^x(t)$ and $\dot{\ox}(t)$ are parallel thanks to \eqref{ort}, i.e.,
\begin{equation}\label{q}
q^x(t)=C\dot{\ox}(t)
\end{equation}
for some constant $C$.

Moreover, it follows from conditions {\bf(5)} and {\bf(7)} that
$$
\lm\oa=8\bigg\la q^x(t),\dfrac{\ox(t)}{\|\ox(t)\|}\bigg\ra,
$$
and thus
$$
\la q^x(t),\dot{\ox}(t)\ra=-8\oa(t)\l\la q^x(t),\dfrac{\ox(t)}{\|\ox(t)\|} \r\ra=-\lm[\oa]^2
$$
since $\left\la q^x(t),\ox^{obs} -\ox(t)\right\ra=0$. This also tells us that $q^x(t)$ point in the opposite direction of the velocity $\dot\ox(t)$, i.e. $C<0$.
Combining it with \eqref{q} gives us $C\|\dot{\ox}(t)\|^2=-\lm[\oa(t)]^2,$ i.e.,
\begin{equation}\label{K}
\|\dot{\ox}(t)\|=K\|\oa\|,
\end{equation}
where $K:=-\dfrac{\lm}{C}$. For simplicity, we may select $\lm =-C$ to get $K=1$. \vspace{1ex}

\begin{figure}[ht]
\centering
\includegraphics[scale=0.38]{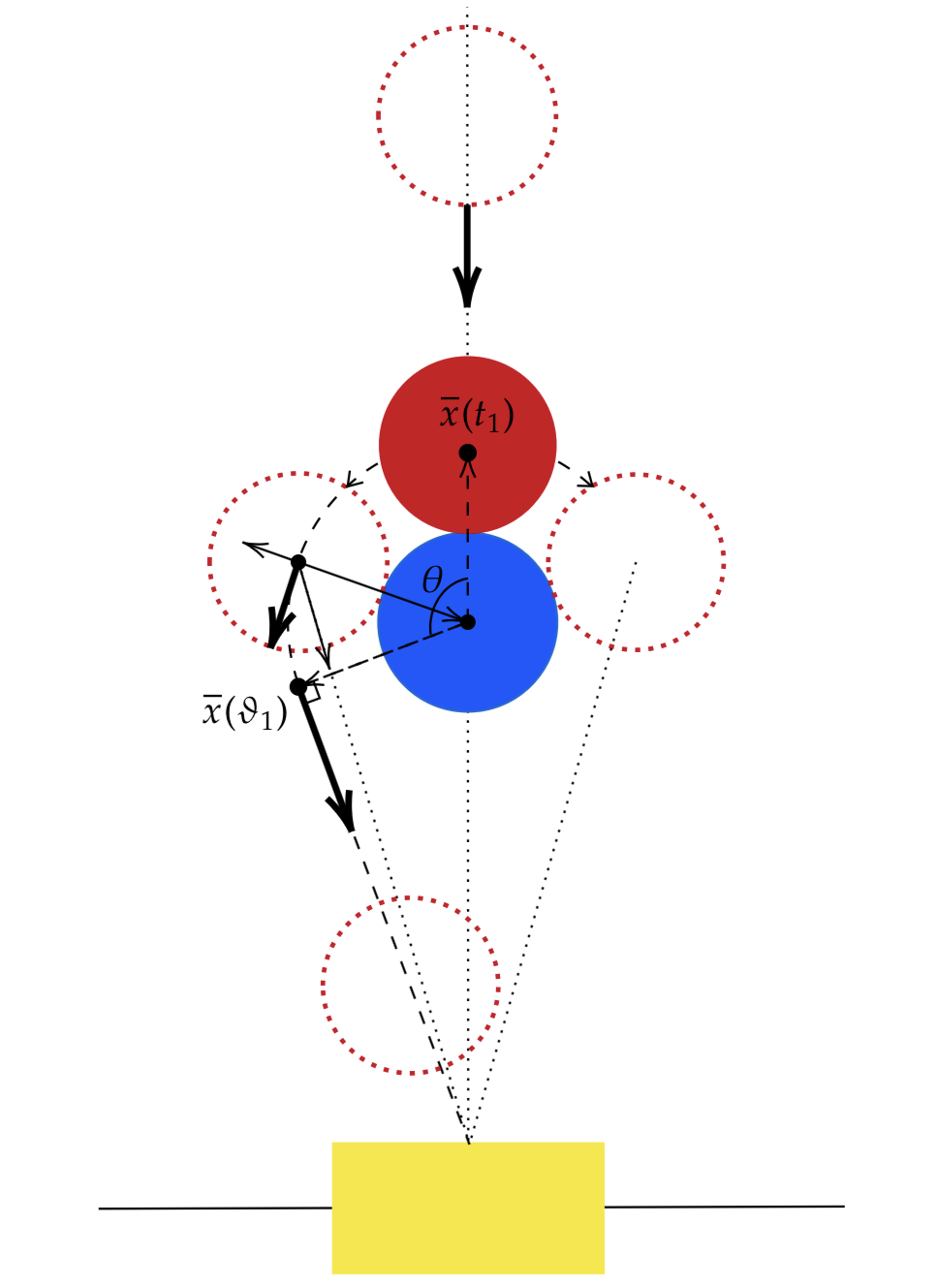}
\caption{The agent is contacting the obstacle.}
\label{pic:general case}
\end{figure}

At the time $t=\vartheta_1$, it makes sense to expect that the agent will leave the obstacle and will head towards the target. In fact, the time $t=\vartheta_1$ can be considered as the safe time for the agent to move to target using his/her own desired velocity. The question is how to track the position of the agent at this very moment, i.e., $\ox(\vartheta_1)$. It seems that the normal cone $N_{C^{obs}}(\ox(\vartheta_1))$ in this case is no longer active, which means that $\eta(\vartheta_1)=0$. Fortunately, equation \eqref{eta} for $\eta(\cdot)$ tells us that $\eta(\vartheta_1)=0$ if
$$
\l\la \ox(\vartheta_1), \ox(\vartheta_1) -x^{obs}\r\ra =0,
$$
which allows us to track the position of the agent at the point where two vectors $\ox(\vartheta_1)$ and $\ox(\vartheta_1) -x^{obs}$ are orthogonal (see Figure~\ref{pic:general case}). This useful observation leads us to valuable information about the time $\vartheta_1$.

To proceed further, let $\theta$ be the angle of two vectors $\ox(t_1)-x^{obs}$ and $\ox(\vartheta_1) - x^{obs}$.
 Then
\begin{align*}
\th&=\cos^{-1}\bigg(\bigg\la\vt j,\,\dfrac{\ox(\vth_1)-\ox^{obs}}{\|\ox(\vth_1)-\ox^{obs}\|} \bigg\ra\bigg) \nonumber\\
&=\cos^{-1}\bigg(\bigg\la (0,1),\,\dfrac{\ox(\vth_1)-(0,24)}{\|\ox(\vth_1)-(0,24)\|} \bigg\ra\bigg),
\end{align*}
where $\vt j = (0,1)$ is the standard unit vector. Observe that $\ox(\vth_1)$ is the intersection of the circle $\mathcal (C)$ centered at $\ox^{obs}$ with the radius $r=\ox(\vth_1)-\ox^{obs}$ and the tangent line to $\mathcal C$ through the target $(0,0)$.

It is not hard to see that there are two possible values for $\ox(\vth_1)$ as follows: either
$$\ox(\vth_1)=\bigg(-\dfrac{3}{2}\sqrt{15},\,\dfrac{45}{2}\bigg),\;\mbox{ or}
$$
 $$
 \ox(\vth_1)=\bigg(\dfrac{3}{2}\sqrt{15},\,\dfrac{45}{2}\bigg).
 $$
As a consequence, we have
$$
\th=\cos^{-1}\bigg(\bigg\la (0,1),\,\bigg(\pm\dfrac{\sqrt{15}}{4},-\dfrac{1}{4}\bigg)\bigg\ra\bigg)\approx 104.48^{\circ}.
$$
Thus the distance of agent's travel from $t=t_1$ to $t=\vth_1$ is
$$
\int_{t_1}^{\vth_1}\|\dot{\ox}(t)\|dt=\dfrac{\pi 2R\th}{180^\circ}.
$$
Combining this with \eqref{K}, we come up to
$$
\int_{t_1}^{\vth_1}\oa dt \approx 10.9411,
$$
and hence obtain
$$
\vth_1=\frac{10.9411}{\oa}+t_1.
$$
In addition, we can determine $t_1$ by
$$
t_1=\frac{24}{8|\oa |},
$$
which leads us in turn to
$$
\vth_1=\frac{13.9411}{|\oa |}.
$$
Meanwhile, it is not hard to observe that $\l\| \ox(\vth_1)\r\| =\sqrt{540}$ (see Figure~\ref{pic:general case}). The distance that the agent travels from $t=\vth_1$ to $t=6$ (the ending time) is $\l\| \ox(\vth_1)\r\| - \| \ox(6)\|$, which can also be found by multiplying the travel time $(6-\vth_1)$ by the constant controlled speed $8\oa$.
As a result, we
$$
(6-\vth_1)8\oa=\sqrt{540}-\|\ox(6)\|,
$$
and get therefore that
$$
\|\ox(6)\|=\sqrt{540}+(\vth_1-6)8\oa,
$$
which gives us the following computation of the cost functional:
\begin{align*}
J&=\dfrac{1}{2}\big[\sqrt{540}+(\vth_1-6)8\oa\big]^2+3\oa^2\nonumber\\
&=\dfrac{1}{2}\bigg[\sqrt{540}+\bigg(\dfrac{13.9411}{\oa}-6\bigg)8\oa\bigg]^2+3\oa^2.
\end{align*}
This quadratic cost function achieves its minimum value $J\approx 23.5871$ at $\oa\approx 2.8004$ with the contacting times $t_1\approx 1.0713$ and $\vth_1\approx 4.9783$.

Furthermore, the trajectory $x(\cdot)$ can be computed by integration as follows:
$$\ox_1(t)\approx\nn
\begin{array}{lll}
&\left(0,48-16.8018t\right),\;\textrm{if}\;t\in (0,t_1),\\
&\left(-6\sin\theta(t),6\cos \theta(t)+24\right),\\
&\textrm{where}\;\theta(t)\in(0,104.48^{\circ}),\;\textrm{if}\;t\in [t_1,\vth_1],\\
&\left(-10.7878+t,41.7809-\sqrt{15}t\right),\\
&\textrm{if}\;t\in (\vth_1,6).
\end{array}
\right.$$
Note that adjusting the value of $\tau>0$ in the cost functional will affect the performance of the agent. This is shown in the table below.

\begin{center}
    \begin{tabular}{|l|c|c|c|c|}
        \hline
        $\tau$ & $\oa$ & $t_1$ & $\vth_1$ & $J[\ox, \ou]$\\
        \hline
        1.0 & $ 2.8004   $ & $1.0713  $ & $4.9783$ & $23.5871$ \\
        1.5 & $ 2.7967   $ & $ 1.0727  $ & $4.9848$   & $35.3348$ \\
        2.0 & $  2.7931  $ & $1.0741   $ & $4.9913$ &$47.052$ \\
        2.5 & $2.7895 $ & $ 1.0755 $ & $4.9977$ & $ 58.7389$ \\
        3.0 & $ 2.7859 $ & $ 1.0769$ & $5.0042$ & $70.3956$ \\
        3.5 & $ 2.7823    $ & $ 1.0783 $ & $ 5.0107$ & $82.0222 $ \\
        4.0 & $  2.7787   $ & $1.0797  $ & $5.0172 $ & $93.6189 $ \\
        4.5 & $ 2.7751    $ & $ 1.0810 $ & $ 1.0810$ & $105.1857 $ \\
        5.0 & $ 2.7716    $ & $ 1.0824 $ & $ 5.0301$ & $ 116.7228$ \\
        5.5 & $ 2.7680    $ & $ 1.0838 $ & $ 5.0365$ & $ 128.2302$ \\
        6.0 & $ 2.7645    $ & $ 1.0852 $ & $ 5.0429$ & $ 139.708$ \\
        6.5 & $ 2.7609    $ & $ 1.0866 $ & $ 5.0495$ & $ 151.157$ \\
        7.0 & $ 2.7574    $ & $ 1.0880 $ & $ 5.0559$ & $ 162.576$ \\
        7.5 & $ 2.7539    $ & $ 1.0894 $ & $ 5.0623$ & $ 173.966$ \\
        8.0 & $ 2.7503    $ & $ 1.0908 $ & $ 5.0689$ & $ 185.327$ \\
        8.5 & $ 2.7468    $ & $ 1.0922 $ & $ 5.0754$ & $ 196.659$ \\
        9.0 & $ 2.7433    $ & $ 1.0936 $ & $ 5.0819$ & $ 207.963$ \\
        9.5 & $ 2.7399    $ & $ 1.0949 $ & $ 5.0882$ & $ 219.237$ \\
        10.0 & $ 2.7364   $ & $ 1.0963 $ & $ 5.0947$ & $ 230.483$ \\
        \hline
    \end{tabular}
\end{center}\vspace*{0.1in}

It occurs  therefore that increasing the value of $\tau$ makes the agent use less energy, but this also means that he/she would stay in contact with the obstacle longer. This example clearly shows that the crowd motion models with obstacles are significantly more challenging  in comparison with the models without obstacles considered in \cite{cm1}.

 \section{Future Research}
 In our future research, we intend to study a more general and more realistic framework of crowd motion models with multiple obstacles, develop necessary optimality conditions and numerical algorithms to solve the corresponding optimal control problems for the obtained sweeping dynamics.

\section{Acknowledgment}
The author Tan H. Cao acknowledge the support of the National Research Foundation of Korea grant funded by the Korea Government (MIST) NRF-2020R1F1A1A01071015.


The author Boris S. Mordukhovich acknowledge the support of the USA National Science Foundation under grants DMS-1007132 and DMS-1512846, by the USA Air Force Office of Scientific Research grant \#15RT0462, and by the Australian Research Council under Discovery Project DP-190100555.

The authors Dao Nguyen and Trang Nguyen acknowledge the support of the USA National Science Foundation under grant DMS-1512846 and by the USA Air Force Office of Scientific Research grant \#15RT0462.
\bibliographystyle{apalike}

\end{document}